\documentclass[11pt]{amsart}
\usepackage{amsbsy,amsfonts}
\usepackage{latexsym,eucal,amsmath,amsthm}
\usepackage{amssymb}

\newcommand{\half}{\frac{1}{2}}

\newcommand{\lam}{\lambda}

\newcommand{\R}{{\mathbb R}}

\newcommand{\C}{{\mathbb C}}

\newcommand{\Z}{{\mathbb Z}}
\newcommand{\N}{{\mathbb N}}

\newcommand{\Schwarz}{\mathcal S}

\newcommand{\Lp}[2]{{\left\| #2 \right\|}_{L^ #1 }}

\newcommand{\LTwo}[1]{{\left\| #1 \right\|}_{L^2}}

\newcommand{\LTwoxT}[1]{{\left\| #1 \right\|}_{L^2_{xT}}}

\newcommand{\Hsub}[2]{{{\left\| #2 \right\|}_{{H_#1}}}}
\newcommand{\Hsup}[2]{{{\left\| #2 \right\|}_{{H^#1}}}}

\newcommand{\X}[3]{{{\left\| #3 \right\|}_{{X_{#1,#2}}}}}

\newcommand{\Xsb}{{X_{s,b}}}

\newcommand{\Boxinv}{ {\Box}^{-1} }
\newcommand{\BID}{ \Boxinv \Delta }
\newcommand{\BIDhalf}{ {{ (\BID) }^\half} }

\newcommand{\Ksub}[1]{{(1+ |{k_ #1 }|)}}

\newcommand{\q}{{|\lambda +{{|k|}^2}|}}

\newcommand{\qsubpm}[1]{{|{\lambda_ #1} \pm {{|{k_ #1 }|}^2}|}}

\newcommand{\w}{{|\lam \pm |k||}}

\newcommand{\ctermpm}[3]{\frac {{ { {{\Ksub #1 }^{ #2 }} c ({k_
                        #1},{\lam_ #1 }) }}}
                        {{(1+|{\lambda_ #1 } \pm {{|{k_ #1 }|}^2}|) }^{ #3 }}}

\newcommand\Conv{  {{
                        {\begin{matrix}
                           {\scriptscriptstyle{{k={k_1}+{k_2}}}}  &\\
                           {\scriptscriptstyle{{\lam={\lam_1}+{\lam_2} }}} &
                         \end{matrix}  }}}  }

\newcommand\siglims{  {{
                        {\begin{matrix}
                 {\scriptscriptstyle{0 \leq \sigma_1 , \sigma_2, \sigma_3 
                            \leq s  }} & \\
                 {\scriptscriptstyle{ \sigma_1 + \sigma_2 + \sigma_3 = s}} &
                         \end{matrix}  }  }}  }

\newtheorem{theorem}{Theorem}
\theoremstyle{definition}

\theoremstyle{remark}
\newtheorem{remark}{Remark}
\theoremstyle{proposition}
\newtheorem{proposition}{Proposition}
\theoremstyle{lemma}
\newtheorem{lemma}{Lemma}
\theoremstyle{corollary}

\numberwithin{equation}{section}
\numberwithin{proposition}{section}
\numberwithin{definition}{section}
\numberwithin{corollary}{section}

\begin{document}

\title{Regularity Bounds on Zakharov System Evolutions}
\author{J.~Colliander}
\thanks{J.C. was supported in part by N.S.F. Grant DMS 0100595.}
\address{\small University of Toronto}
\author{G.~Staffilani}
\thanks{G.S. was supported in part by N.S.F. Grant 
DMS 0100375 and
grants from Hewlett and Packard and the Sloan Foundation.}
\address{\small Brown University and Stanford University}

\begin{abstract}
Spatial regularity properties of certain global-in-time
solutions of the Zakharov system are established. In particular, the
evolving solution $u(t)$ is shown to satisfy an estimate $\Hsup s {u(t)}
\leq C {{|t|}^{(s-1)+}}$, where $H^s$ is the standard spatial Sobolev
norm. The proof is an adaptation of earlier work on the nonlinear
Schr\"odinger equation which reduces matters to bilinear estimates.
\end{abstract}

\date{14 March 2002}

\maketitle

\section{introduction}

We consider the initial value problem for the Zakharov system 
\cite{Z} on $\R^2$
\begin{equation}
\left\{ \begin{matrix}
i{u_t} + \Delta u =nu, & u : \R^2 \times [-T_*, T^*] \longmapsto \C , \\
\Box n = \Delta  {{|u|}^{2}}, 
& n: \R^2 \times [ - T_*, T^*]\longmapsto \R , \\
(u , n , \dot{n} )(0) = (\phi, a , b).
\end{matrix} \right.
\label{ZS} 
\end{equation}
Suppose $b$ is such that there exists $V: \R^2 \rightarrow \R^2$ with
$b = \nabla \cdot V$. Then the Zakharov system may be rewritten in Hamiltonian
form with Hamiltonian
\begin{equation}
H(u, \bar{u} , n , V) = \int_{\R^2} {{|\nabla u |}^2}
+ \half ( n^2 + {{|V|}^2}) + n {{|u|}^2} dx.
\label{Hamiltonian}
\end{equation}
For initial data $\phi$ small enough in $L^2$ we can condlude from conservation
of \eqref{Hamiltonian} that
\begin{equation}
\Hsub 1 {{u(t), n(t) , \dot{n} (t) }} =
{{\| u(t) , n(t) , \dot{n}(t) \|}_{H^1 \times L^2 \times {{\widehat{H}}^{-1}} }} \leq C \Hsub 1 {\phi, a , b }
\label{aprioriHone}
\end{equation}
where $H_1 := H^1 \times L^2 \times {{\widehat{H}}^{-1}}$ and $
{{\widehat{H}}^{-1}}$ is defined by ${{\| b \|}_{ {{\widehat{H}}^{-1}} }} =
\Lp 2 V $.
Local wellposedness of \eqref{ZS} for data $(\phi, a , b ) \subset H_1$ was
established in \cite{BouCol97}, \cite{Thesis}, with the lifetime of existence
satisfying
\begin{equation}
T > {{ \|{\phi, a , b }\| }_{H_1}^{-\alpha}} ~{\mbox{for some}}~ \alpha > 0.
\label{lifetime}
\end{equation}
(The regularity requirements for the local results in \cite{BouCol97}, 
\cite{Thesis}
have subsequently been improved in \cite{GTV}.) Hence, for data implying
a priori $H_1$ control \eqref{aprioriHone}, the local result may be iterated
to prove the existence of global-in-time solutions of \eqref{ZS}. In fact,
global solutions of the initial value problem \eqref{ZS} had been shown to 
exist earlier  \cite{AA} ($d=1$)  \cite{SS} ($d=2$) using energy methods in spaces requiring more regularity than $H_1$.

\begin{remark} The initial value problem \eqref{ZS} has solutions which
blow up in finite time \cite{GM1}, \cite{GM2}. At present, there is no
criteria known which identifies those initial data leading to blow up
and those leading to global-in-time solutions. This paper provides
regularity bounds on those global solutions obtained by iterating the 
local wellposedness result using a priori $H_1$ control.
\end{remark}

Let $\Boxinv F$ denote the solution of the inhomogeneous wave equation
with zero data,
\begin{equation}
\left\{ \begin{matrix}
\Box n = F, \\
(n , \dot{n} )(0) = (0 , 0).
\end{matrix} \right.
\label{inhomog} 
\end{equation}
Let $W(a,b)$ denote the solution of
\begin{equation}
\left\{ \begin{matrix}
\Box n = 0 \\
(n , \dot{n} )(0) = (a , b).
\end{matrix} \right.
\label{homog} 
\end{equation}
Note that $\Boxinv F$ and $W(a,b)$ may be explicitly represented using the 
Fourier transform. The (formal) solution of the second equation in \eqref{ZS}
is
\begin{equation}
n = W(a,b) + \Boxinv ({ \Delta {{|u|}^2} }).
\label{nsolved}
\end{equation}
Substituting this expression for $n$ into the first equation in \eqref{ZS}
gives
\begin{equation}
\left\{ \begin{matrix}
u_t = i \Delta u - i W(a,b) u - i \Boxinv \Delta ( {{|u|}^2} ) u, \\
u (0) = \phi .
\end{matrix} \right.
\label{nremoved} 
\end{equation}
Note that the regularity properties of the data $a,b$ and of $u$, 
inferred from solving 
\eqref{nremoved}, determine the regularity properties of $n$ through
\eqref{nsolved}.

Let $\Schwarz$ denote the Schwarz class. Consider initial data $\phi, a \in \Schwarz, ~ b \in \Schwarz \bigcap
{{\widehat{H}}^{-1}}$ implying a priori $H_1$ control \eqref{aprioriHone}.
How do the regularity properties of the global solution $(u(t), n(t))$
behave as $t \rightarrow \infty$? In particular, can we describe, or
at least bound from above and below, $\Hsub s {{u(t) , n(t) , \dot{n} (t) }} $
for $s \gg 1$ as $t \rightarrow \infty$? These estimates quantify the shift of the conserved $L^2$ mass in frequency space. In particular, the upper bounds we obtain in this paper limit the rate of transfer from low frequencies to high frequencies. By the note following
\eqref{nremoved}, it suffices to understand $\Hsup s {u(t) }$.

The local result for \eqref{nremoved} implies
\begin{equation}
\sup_{t \in [0,T]} \Hsup s {u(t) } \leq \Hsup s {\phi} + C \Hsup s {\phi}
\label{givesexp}
\end{equation}
which iterates to give an exponential bound $\Hsup s {u(t) } \leq C^{|t|}.$
Bourgain observed that a slight improvement of \eqref{givesexp}
\begin{equation}
\sup_{t \in [0,T]} \Hsup s {u(t) } \leq \Hsup s {\phi} + C {{ \| \phi \|}_{
H^s}^{1 - \delta}}, ~0 < \delta < 1,
\label{givespoly}
\end{equation}
implies the polynomial bound 
$\Hsup s {u(t) } \leq C {{|t|}^{\frac{1}{\delta}}}.$ This observation was
exploited in \cite{B4} to prove polynomial bounds on high Sobolev norms
for solutions of the nonlinear Schr\"odinger equation and certain 
nonlinear wave equations.

Staffilani \cite{Staff97}, \cite{StaffQ97} improved the degree of the 
polynomial upper bound using a different approach to prove \eqref{givespoly}
in the case of the nonlinear Schr\"odinger equation. The crucial bilinear 
estimate used in this approach has recently been improved \cite{CDKS}
giving a slightly better polynomial estimate. This paper adapts the arguments
from \cite{StaffQ97} for the nonlinear Schr\"odinger equation to 
prove similar polynomial bounds on high Sobolev norms for the global
solutions of \eqref{ZS} constructed in \cite{BouCol97}, \cite{Thesis}.

\begin{theorem}
\label{zspoly}
Assume $(\phi, a ,b ) \in \Schwarz \times \Schwarz \times ( \Schwarz \bigcap
{{\widehat{H}}^{-1}})$. 
Global solutions of \eqref{ZS} satisfying \eqref{aprioriHone} also satisfy
\begin{equation}
\Hsup s {u(t)} \leq C {{|t|}^{(s-1)+}} .
\end{equation}

\end{theorem}

The question of lower bounds showing growth of high Sobolev norms 
remains a fascinating open question. For a more thorough discussion, including
model equations other than the Zakharov system, see the book of
Bourgain \cite{BAMS}.

\section{Reduction to bilinear estimate}

Our goal is to bound $\Hsup s {u(t)}$ for $u$, the solution of 
\eqref{nremoved}, with $t \in [0,T]$ and $T$ as in \eqref{lifetime}.
Since $\Lp 2 {u(t)} = \Lp 2 {\phi}$ for all $t$, it suffices to bound
$\Lp 2 {B^s u(t) }$ where $B = {\sqrt{-\Delta}}$. Let's assume $s = 2 m$
for $1 \ll m \in \Z$ to avoid certain technical issues involving fractional
derivatives below. Let $\langle \cdot, \cdot \rangle$ denote the 
standard $L^2_x$
inner product, $\langle f , g \rangle = \int_{\R^2}  f {\bar{g}} dx.$
By the fundamental theorem of calculus,
\begin{equation}
{{\| B^s u(t) \|}_{L^2}^2} = {{\| B^s u(0) \|}_{L^2}^2} + 
\int_0^t \frac{d}{d\sigma} \langle B^s u(\sigma ) , B^s u(\sigma) \rangle 
d\sigma .
\label{one}
\end{equation}
We calculate
\begin{equation}
I = 2 \Re \int_0^t \langle B^s \dot{u} (\sigma ) , B^s u(\sigma ) \rangle
d \sigma .
\label{two}
\end{equation}
Now, using the equation \eqref{nremoved}, we find
\begin{equation}
I = - 2 \Im \int_0^t \langle B^s \Delta u (\sigma ), B^s u(\sigma ) \rangle
d \sigma +
\label{three}
\end{equation}
\begin{equation*}
2 \Im \int_0^t \langle B^s [ W(a,b)(\sigma) u(\sigma)] , 
B^s u(\sigma) \rangle d \sigma + 2 \Im \int_0^t \langle B^s [ \Boxinv \Delta
({{|u|}^2}) u (\sigma )], B^s u(\sigma ) \rangle d \sigma .
\end{equation*}
Denote the three terms on the right-side of \eqref{three} by $I_1 + I_2 + I_3$.
Upon writing $- \Delta = B^2$ and integrating by parts, the first term $I_1$
is seen to have a real integrand so this term is zero. The term $I_2$
involves $B^s (W(a,b) u )$. Various terms arise from the Leibniz rule
for differentiating a product. The most dangerous of these is $W(a,b) B^s u$
but, since $W(a,b) $ is a real-valued function, this term leads to a purely
real integrand in \eqref{three} and so disappears. Hence, the term $I_2$
leads to a sum of terms of the form
\begin{equation}
C \Im \int_0^t \langle [B^{s_1} W(a,b)(\sigma)][B^{s_2} u(\sigma )] , 
B^s u(\sigma ) \rangle d \sigma,
\label{four}
\end{equation}
where $s = s_1 + s_2 , ~ 1 \leq s_1 \leq s, ~0 \leq s_2 \leq s-1, ~
s_1, s_2 \in \Z$,

We can multiply by a smooth cutoff function in time $\psi_T \thicksim
\chi_{[0,T]}$ and estimate these terms via the H\"older inequality by
\begin{equation}
{{\| B^{s_1} W(a,b) \|}_{L^2_{x,t \in [0,T]}}} {{\| B^{s_2} u \|}_{L^4_{xt}}}
{{\| B^s u \|}_{L^4_{xt}}}.
\label{five}
\end{equation}
The Strichartz estimate for the paraboloid and properties of $\Xsb$ spaces{\footnote{The space $X_{s,b}$ is defined using the norm
$  {{\| u \|}_{X_{s,b}}} = {{ \left( \int (1+|k|)^{2s} (1+|\lam + |k|^2|)^{2b}
| \widehat{u} (k , \lam )|^2 dk d\lam \right)}^\half}.$}}
\cite{B1} imply for $ b = \half+$,
\begin{equation}
{{\| B^{\tilde{s}} u \|}_{L^4_{xt}}} \leq C_T \X {\tilde{s}} b u .
\label{six}
\end{equation}
% The space $X_{s,b}$ is defined using the norm
% \begin{equation}
%   \label{Xsbnorm}
%   {{\| u \|}_{X_{s,b}}} = {{ \left( \int (1+|k|)^{2s} (1+|\lam + |k|^2|)^{2b}
% | \widehat{u} (k , \lam )|^2 dk d\lam \right)}^\half}.
% \end{equation}
The local wellposedness result \cite{BouCol97}, \cite{GTV} gives
\begin{equation}
\X {\tilde{s}} b u 
\leq C \Hsup {{\tilde{s}}} {u(0)}
\label{seven}
\end{equation}
Therefore, the second term $I_2$ in \eqref{three} is estimated by a sum
of terms of the form
\begin{equation}
{{\| a, b \|}_{H^{s_1} \times H^{{s_1} - 1}}} \Hsup {{s_2}} {u(0)}
\Hsup s {u(0)} .
\label{eight}
\end{equation}
The first factor is bounded by a constant which depends upon the initial
data $a,b$. The second factor may be interpolated between $\Hsup 1 {\phi} $
and $\Hsup s {\phi} $ which leads to the bound
\begin{equation}
|I_2 | \leq C \sum\limits_{0 \leq s_2 \leq s-1 } 
{{\| u(0) \|}_{H^s}^{1 + {\frac{{s_2} -1}{s-1}}}} \leq C 
{{\| u(0) \|}_{H^s}^{2 - {\frac{1}{s-1}}}}.
\label{nine}
\end{equation}

It remains to bound $I_3$. Since differentiation in $x$ commutes with
$\BID$, the Leibniz rule shows
\begin{equation}
I_3 = 2 \Im \sum\limits_\siglims
c_{\sigma_1 , \sigma_2 , \sigma_3 }
\int_0^t \langle \BID ( B^{\sigma_1} u B^{\sigma_2} {\bar{u}} ) 
B^{\sigma_3} u , 
B^s u \rangle d \tau.
\label{ten}
\end{equation}
In case $\sigma_3 = s$, the resulting integrand is purely real so this
term disappears. Consider first those terms with $\sigma_3 \leq s - 2$
and after treating these we will consider the terms with $\sigma_3 = s-1$.
Since we are interested in proving a local-in-time estimate, we can
insert a smooth cutoff $\psi_T \thicksim \chi_{[0,T]}$ and wish to 
bound
\begin{equation}
\label{eleven}
| \psi_T (t) \int_0^t \BID ( B^{\sigma_1} u B^{\sigma_2} {\bar{u}} )
B^{\sigma_3} u B^s {\bar{u}} d\tau |.
\end{equation}
A formal ``integration by parts'' (which is justified rigorously in the
next section when we define $\BIDhalf$) allows us to bound by
\begin{equation}
\label{twelve}
| \psi_T (t) \int_0^t \BIDhalf ( B^{\sigma_1} u B^{\sigma_2} {\bar{u}} )
\BIDhalf( B^{\sigma_3} u B^s {\bar{u}}) d\tau |
\end{equation}
and Cauchy-Schwarz reduces matters to controlling
\begin{equation}
\LTwoxT {  \BIDhalf ( B^{\sigma_1} u B^{\sigma_2} {\bar{u}} ) }
\LTwoxT {  \BIDhalf( B^{\sigma_3} u B^s {\bar{u}})  } .
\label{thirteen}
\end{equation}

\begin{proposition}
\label{bilinearestimate}
Let $0 \leq s_1 \in \N$ and $1 \leq s_2 \in \N,~ s_1 \leq s_2.$ For
$b > \half$,
\begin{equation}
\label{propstatement}
\LTwoxT { \BIDhalf ([B^{s_1} u_1 ][ B^{s_2} {\bar{u_2}} ]) }
\leq C \X {{s_1}+1} b {u_1}  \X {s_2 - \half } b {u_2} .
\end{equation}
The estimate is also valid if the complex conjugation is moved from
$u_2$ to $u_1$ on the left-side of \eqref{propstatement}.
\end{proposition}

Suppose the proposition is true. The bilinear expressions in \eqref{thirteen}
are estimated
\begin{equation}
\X {\sigma_1 + 1} b {u}  \X {\sigma_2 - \half} b u
\X {\sigma_3 + 1 } b u  \X {s - \half } b u .
\label{fifteen}
\end{equation}
Using the local result we know $ \X {\tilde{s}} b u \leq 
C \Hsup {{\tilde{s}}} {u(0)} $ 
and upon interpolating the various $H^{\tilde{s}}$
norms between $H^1$ and $H^s$ (using \eqref{aprioriHone}) bounds \eqref{fifteen}
by 
\begin{equation}
C {{\| u(0) \|}_{H^s}^{ \frac{\sigma_1 + \sigma_2 - \half - 1 + \sigma_3 + s 
- \half - 1 }{s-1}}} .
\label{sixteen}
\end{equation}
Recalling that $\sigma_1 + \sigma_2 + \sigma_3 = s,$ the exponent simplifies
to $2 - \frac{1}{s-1}$, just as in \eqref{nine}.

Now, consider a term in \eqref{ten} with $\sigma_3 = s-1 $. Evidently,
$\sigma_1 = 1,~ \sigma_2 = 0$ or $\sigma_1 = 0,~\sigma_2 = 1$. In this
case, we apply Cauchy-Schwarz directly to the term as it appears in 
\eqref{ten} to bound by
\begin{equation}
\label{seventeen}
\LTwoxT { \BID ( B^{\sigma_1 } u B^{\sigma_2} {\bar{u}} ) }
\LTwoxT { B^{\sigma_3} u B^s {\bar{u}} } .
\end{equation}
The second factor is readily estimated using Bourgain's refinement of the
Strichartz inequality \cite{BRefine} to give
\begin{equation}
\label{eighteen}
\X {\sigma_3 + \half +} b u  \X {s - \half } b u \leq 
C {{\| u \|}_{{X_{s-\half + , b}}}^2}.
\end{equation}
The first factor in \eqref{seventeen} is bounded using a variant of
Proposition \ref{bilinearestimate}.

\begin{proposition}
\label{nregularity}
Let $0 \leq s_1 \in \N,~1 \leq s_2 \in \N,~ s_1 \leq s_2$. For $b > \half$,
\begin{equation}
\label{nineteen}
\LTwoxT { \BID ( [B^{s_1} u_1 ] [ B^{s_2} {\bar{u_2}} ] ) }
\leq C \X {s_1 + 1 } b {u_1 }  \X {s_2} b u .
\end{equation}
The estimate is also valid if the complex conjugation is moved from
$u_2$ to $u_1$ on the left-side of \eqref{nineteen}.
\end{proposition}

Combining \eqref{eighteen}, \eqref{nineteen} shows \eqref{seventeen}
may be bounded using $\sigma_1 = 1, \sigma_2 = 0$ or $\sigma_1 = 0, \sigma_2 
= 1$ and $\sigma_3 = s-1$,
\begin{equation*}
\X 1 b u  \X 1 b u  {{\| u \|}_{{X_{s-\half + , b}}}^2}
\leq C  {{\| u \|}_{{X_{s-\half + , b}}}^2} .
\end{equation*}
The local result and interpolation bounds this by
\begin{equation*}
C{{\| u(0) \|}_{H^s}^{ 2 ( \frac{s - \half - 1 }{s-1} )+}} = 
C{{\| u(0) \|}_{H^s}^{ 2- \frac{1}{s-1} +}}
\end{equation*}
which (up to the $+$) is the same as in \eqref{nine}, \eqref{sixteen}.

Summarizing, the two Propositions above show that the integral term in 
\eqref{one} is bounded by
\begin{equation*}
C{{\| u(0) \|}_{H^s}^{ 2- \frac{1}{s-1} +}}.
\end{equation*}
We may assume that $\Lp 2 { {B^s} u(t) } \geq \Lp 2 { {B^s} u(0) }$ for
otherwise $\eqref{givespoly}$ is automatic. Therefore, we can divide 
\eqref{one} through by $\Lp 2 {B^s u(t)} $ and with $L^2$ conservation
observe \eqref{givespoly} holds with $\delta = \frac{1}{s-1}-$ 
proving Theorem \ref{zspoly}.

The next section establishes the Propositions and defines $\BIDhalf$ used 
in the treatment of $\sigma_3 \leq s-2$ terms above.

\section{Bilinear Estimates}

In this section, we present a proof of Proposition \ref{bilinearestimate}.
Along the way we will observe explicit properties of the operator
$\BID$ which allow us to justify step \eqref{twelve} in the previous
section. Proposition \ref{nregularity} will follow from modifications of
the proof of Proposition \ref{bilinearestimate}.

The operator $\Boxinv$ was defined as the mapping taking the inhomogeneity
$F$ to the solution of the linear initial value problem \eqref{inhomog}.
It can be explicitly represented using the Fourier transform as
\begin{equation}
\Boxinv F (x,t) = - \int \int
e^{i k \cdot x } \left[ e^{i \lambda t } - \half
( 1 + \frac{\lambda}{|k|} ) e^{i |k| t } - \half
( 1 - \frac{\lambda}{|k|} ) e^{-i |k| t }  \right]
\frac{ {\widehat{F}}(k , \lam ) }{ (\lam - |k|) (\lam + |k|) } dk d\lam 
\label{bone}
\end{equation}
where ${\widehat{F}}$ denotes the space-time Fourier transform of $F$.
A Taylor series argument show that the apparent singularities along $\lam
\pm |k| = 0$ do not occur and that
\begin{equation}
\left| {\widehat{\Boxinv \Delta F}} (k, \lam) \right|
\leq C \left| {\widehat{F}} (k , \lam ) \right| \frac{ |k| }{ ( 1 + |\lam \pm
|k| |) } + \left| {\widehat{F}} (k , \lam ) |_{\{\lam = \mp |k|\} }  \right|
\frac{ |k| }{ ( 1 + |\lam \pm |k| |) } .
\label{btwo}
\end{equation}
From an $L^2$ point-of-view, it is therefore natural to define for 
real numbers $\alpha$,
\begin{equation}
{\widehat{ {{(\BID )}^\alpha} } } ( k, \lam ) =
\left| {\widehat{F}} (k , \lam ) \right| {{\left(\frac{ |k| }{ ( 1 + |\lam \pm
|k| |) } \right)}^\alpha} 
+ \left| {\widehat{F}} (k , \lam ) |_{\{\lam = \mp |k| \}}  \right|
{{ \left(\frac{ |k| }{ ( 1 + |\lam \pm |k| |) } \right)}^\alpha} .
\label{bthree}
\end{equation}
In particular, we have defined the operator ${{ \left( \BID \right)}^\half}$
which appears in the statement of Proposition \ref{bilinearestimate}.
For two functions of space-time, $F,G$, which are cutoff to $t \in [0,T]$,
consider the expression (analagous to \eqref{eleven})
\begin{equation}
\int \int {(\BID ) (F )} ~G dx dt 
= \int \int {\widehat{ {{ (\BID ) (F ) }} }}
{\widehat{G}} dk d\lam.
\label{bfour}
\end{equation}
We insert \eqref{btwo} and take the absolute value under the integral
sign. Then, upon writing
$ \frac{|k|}{(1 + |\lam \pm |k|| )} = 
\frac{{{|k|}^\half}}{{{(1 + |\lam \pm |k|| )}^\half}}
\frac{{{|k|}^\half}}{{{(1 + |\lam \pm |k|| )}^\half}}$ we observe
that
\begin{equation}
\left| \int \int {(\BID ) (F )} ~G dx dt \right|
\leq \left| \int \int {{(\BID )}^\half} \widetilde{F} \cdot
{{(\BID )}^\half} \widetilde{G} dx dt\right|.
\label{bfive}
\end{equation}
where $\widetilde{F} (x,t) = \int e^{i(kx + \lam t)} | \widehat{F} (k,\lam)|
dk d\lam$ and $\widetilde{G}$ is similarly defined. For proving $L^2$-type 
estimates, the distinction between $F$ and $\widetilde{F}$ is
unimportant.  In particular, the ``integration by parts'' 
step \eqref{twelve} is validated.

Now that $\BIDhalf $ has been given a precise meaning, we turn our attention
to proving the inequality \eqref{propstatement}

\begin{proof}[Proof of Proposition \ref{bilinearestimate}]
Since ${\widehat{\BID}} (k, \lam ) \thicksim
\frac{|k|}{\lam \pm |k| }$, we see that $\BID$ can be as bad as
one derivative in $x$. Therefore, the number of derivatives on both
the left-side and right-side of \eqref{propstatement} is $s_1 + s_2 + \frac{1}{2}$.
The desired estimate \eqref{propstatement} may be reexpressed using duality
and certain renormalizations as
\begin{equation}
  \label{reexpressed}
  \int_* d(k, \lam) {{\left( \frac{|k|}{(1+|\lam \pm |k||)}\right)}^\half}
 \ctermpm 1 {-1} b  \ctermpm 2 {\half} b \leq \LTwo d \LTwo {c_1} \LTwo {c_2}
\end{equation}
where $\int_*$ is shorthand for $\int\limits_\Conv $ and without loss we
may assume $d, c_1, c_2 \geq 0$. The choices of $\pm$ in \eqref{reexpressed} 
are assumed to be independent
in the following analysis. In fact, this is only the first contribution
arising from the right-side of \eqref{bthree}. The other ``on-light-cone''
piece may be similarly estimated. We analyze \eqref{reexpressed} in cases
depending upon the size of $|k_1|, |k_2|.$

{\bf{Case 1.}} $|k_1|, |k_2| \leq 10$. \newline
We may ignore $|k|^\half (1+|k_1|)^{-1} (1+|k_2|)^{\half}$ and then drop the
(potentially helpful) wave remnant $(1+|\lam \pm |k||)^{\half}$ to bound the 
left-side of \eqref{reexpressed} by
\begin{equation}
  \label{reduced}
  \int_* d(k, \lam) \frac{c_1 ( k_1, \lam_1 ) }{ {{(1+|\lam_1 \pm |k_1|^2|)}^b} }
 \frac{c_2 ( k_2, \lam_2 ) }{ {{(1+|\lam_2 \pm |k_2|^2|)}^b} }.
\end{equation}
Fourier transform properties show this equals $\langle \widehat{ \mathcal{D} },
\widehat{\mathcal{C}_1} * \widehat{\mathcal{C}_2} \rangle = \langle \mathcal{D} ,
\mathcal{C}_1 \mathcal{C}_2 \rangle $ where $\mathcal{D}, \mathcal{C}_1 , \mathcal{C}_2$ are functions of space and time whose Fourier transforms are $d, \frac{c_1 ( k_1, \lam_1 ) }{ {{(1+|\lam_1 \pm |k_1|^2|)}^b} }, 
 \frac{c_2 ( k_2, \lam_2 ) }{ {{(1+|\lam_2 \pm |k_2|^2|)}^b} }$, respectively.
By H\"older, we can estimate by $\| \mathcal{D} \|_{L^2_{xT}} {\| {\mathcal{C}}_1
\|}_{L^4_{xT}} {\| {\mathcal{C}}_2 \|}_{L^4_{xT}} $ and obtain \eqref{reexpressed}
in this case using Plancherel and the Strichartz inequality for the paraboloid
as written in \cite{B1},
\begin{equation}
{{\left\| \int  \frac{a(k, \lam)}{ {{(1+ \q )}^b} } dk d\lam \right\|}_{
L^4 (\R^2_x \times \R_t )}} \leq {{\| a \|}_{L^2_{k\lam}}},~ b > \half,
\label{bstrichartz}
\end{equation}
The standard steps going from\eqref{reduced} through $L^2 L^4 L^4$ to
\eqref{reexpressed} will be omitted from the discussion below.

{\bf{Case 2.}} $|k_1 | \gtrsim |k_2 |,~ |k_1| \gtrsim 10.$ \newline
The case defining conditions imply $|k| \lesssim |k_1|.$ We again ignore the wave
remnant and use $(1+|k_1|)^{-1}$ to cancel away $|k|^\half$ and $(1+|k_2|)^\half$.
We again encounter \eqref{reduced} and complete this case with the $L^2 
L^4 L^4$ H\"older argument using \eqref{bstrichartz}.

{\bf{Case 3.}} $|k_1| \ll |k_2|,~|k_2| \gtrsim 10 \implies |k_2| \thicksim |k|$. \newline
The numerator $(1+|k_1|)^{-1}$ is not helpful in this case so we exploit the
denominators in \eqref{reexpressed} to cancel $|k|^\half$ and $(1+|k_2|)^\half$.
Since $k = k_1 + k_2, ~\lam = \lam_1 + \lam_2$, the triangle inequality implies
\begin{equation}
\max ( \w , \qsubpm 2  , \qsubpm 1 ) \geq | \pm |k|
+ {{|k_2 |}^2} - {{|k_1 |}^2} | \thicksim {{|k_2 |}^2 } \thicksim {{|k|}^2}.
\label{denoms}
\end{equation}
\newline
{\bf{Case 3.A.}} $\w $ is the max in \eqref{denoms}. \newline
We use the large denominator to cancel $|k|^\half$ and $(1+|k_2|)^\half$ and
proceed as with \eqref{reduced}.

{\bf{Case 3.B.}} $\qsubpm 2$ is the max in \eqref{denoms}. \newline
Most of the large denominator is used to cancel away $|k|^\half$ and 
$(1+|k_2|)^\half$ and we need to control
\begin{equation*}
  \int_* \frac{d(k,\lam)}{(1 + \w )^\half} \frac{c_1 (k_1 , \lam_1)}{(1+ 
\qsubpm 1 )^b} \frac{c_2 (k_2 , \lam_2)}{(1+ 
\qsubpm 2 )^{b-\half}} .
\end{equation*}
Since $b> \half$, so that $b - \half > 2\delta > 0$, and $\qsubpm 2 \gtrsim \w$, we can write
\begin{equation}
  \label{stepping}
  \int_*
\frac{ (1+|k|)^{-\delta} d(k , \lam )}{(1+|\lam \pm |k||)^{\half+}}
\ctermpm 1 {-\delta} b  c_2 (k_2 , \lam_2 ).
\end{equation}
Let $\mathcal{D} (k , \lam) = \frac{ (1+|k|)^{-\delta} d(k , \lam )}{(1+|\lam \pm |k||)^{\half+}}$ and $\mathcal{C} = \ctermpm 1 {-\delta} b$. Then \eqref{stepping}
may be expressed as $\langle \mathcal{D} * \mathcal{C} , c_2 \rangle$ and Cauchy-Schwarz reduces matters to controlling $\| \mathcal{D} * \mathcal{C} \|_{L^2_{xT}}$.
This is accomplished in the following lemma. 

\begin{lemma}
\label{waveschrodingerlemma}
For $b=\half+$ and a fixed small $\delta > 0$
\begin{equation}
\label{waveschro}
\int_{*, |k_2|, |k_1| \geq 10} f(k, \lam) \frac{ (1+|k_2|)^{-\delta} d(k_2 , \lam_2 )}{(1+|\lam_2 \pm |k_2||)^{\half+}} \ctermpm 1 {-\delta} b 
\leq \Lp 2 f \Lp 2 d \Lp 2 {c_1}.
\end{equation}
\end{lemma}
\begin{proof}
Since $\half+$ and $b$ exceed $\half$, the estimate \eqref{waveschro} may be reduced to the ``on-curve'' setting using parabolic (for $c_1$) and light-cone (for $d$)
level set decompositions (see, for example, \cite{CDKS}.) This reduces considerations to showing that
\begin{equation}
  \label{oncurve}
  \int_{|k_2| \geq |k_1| \geq 1} f(x_1 + k_2 , \pm |k_1|^2 \pm |k_2| ) |k_2|^{-\delta}
\psi (k_2) |k_1|^{-\delta} \phi (k_1) dk_1 dk_2 
\leq {{\| f \|}_{L^2_{k\lam}}} {{\| \psi \|}_{L^2_k}} {{\| \phi \|}_{L^2_k}} .
\end{equation}
Consider the piece of the integration on the left-side of \eqref{oncurve} arising from $\{ k_2 : |k_2| \thicksim K_2 (dyadic)\} \times \{ k_1 : |k_1| \thicksim K_1 (dyadic)\}$. We make a change of variables{\footnote{Superscripts refer to vector components.}}
$u^1 = k_1^1 + k_2^1,~u^2 = k_1^2 + k_2^2, ~ v = \pm |k_1|^2 \pm |k_2|$ and we 
assume that the component $k_1^2$ of $k_1$ satisfies $k_1^2  \thicksim |k_1| \thicksim K_1.$ (This may be accomplished by cutting in pie slices and making a rotation 
of coordinates if necessary.) This change of variables followed by Cauchy-Schwarz
shows \eqref{oncurve} is bounded by 
\begin{equation}
\label{cov}
  {{\| f \|}_{L^2_{k\lam}}} K_1^{-\delta} K_2^{-\delta}
\int_{|k_1^1| \leq K_1} \left( \int_{|k_i| \thicksim K_i} 
|\psi(k_2)|^2 |\phi(k_1)|^2 \frac{1}{|J|} dk_1^2 dk_2^1 dk_2^2 \right)^\half 
dk_1^1
\end{equation}
where the Jacobian is
\begin{equation}
  |J| = |2k_1^2 \pm 1| \thicksim K_1.
\end{equation}
We apply Cauchy-Schwarz in $k_1^1$ and pick up an extra factor of $K_1^\half$
which is cancelled away by the Jacobian factor. The small prefactors $K_i^{-\delta}$
allow us to sum over large dyadic scales thereby proving \eqref{oncurve} and 
the lemma.
\end{proof}
The lemma shows that \eqref{stepping} is bounded as claimed in \eqref{reexpressed}
which completes the proof of Proposition \ref{bilinearestimate}.
\end{proof}

{\bf{Case 3.C.}} $\qsubpm 1 $ is the max in \eqref{denoms}. \newline
This case follows with a modification of the argument for Case 3.B.

The proof of Proposition \ref{nregularity} follows the same case structure
as the proof of Proposition \ref{bilinearestimate}. The only difference is 
in the accounting of the extra $\half$ derivative in both sides of \eqref{nineteen}
in comparison with \eqref{propstatement}.

\enddocument